# FACTORIZATION OF POINT CONFIGURATIONS, CYCLIC COVERS, AND CONFORMAL BLOCKS

MICHELE BOLOGNESI AND NOAH GIANSIRACUSA

Abstract. We describe a relation between the invariants of $n$ ordered points in $\mathbb{P}^d$ and of points contained in a union of linear subspaces $\mathbb{P}^{d_1} \cup \mathbb{P}^{d_2} \subseteq \mathbb{P}^d$. This yields an attaching map for GIT quotients parameterizing point configurations in these spaces, and we show that it respects the Segre product of the natural GIT polarizations. Associated to a configuration supported on a rational normal curve is a cyclic cover, and we show that if the branch points are weighted by the GIT linearization and the rational normal curve degenerates, then the admissible covers limit is a cyclic cover with weights as in this attaching map. We find that both GIT polarizations and the Hodge class for families of cyclic covers yield line bundles on $\overline{M}_{0,n}$ with functorial restriction to the boundary. We introduce a notion of divisorial factorization, abstracting an axiom from rational conformal field theory, to encode this property and show that it determines the isomorphism class of these line bundles. Consequently, we obtain a unified, geometric proof of two recent results on conformal block bundles, one by Fedorchuk and one by Gibney and the second author.

## 1. Introduction

In this paper we study configuration spaces of points on the line and in higher-dimensional projective spaces, relating classical constructions to modern results.

### 1.1. Geometric invariant theory.

A fundamental object in classical invariant theory is the ring of invariants for $n$ ordered points in projective space up to projectivity. This is the homogeneous coordinate ring for the natural projective embedding of the GIT quotient $(\mathbb{P}^d)^n /\!\!/_{\underline{c}} \operatorname{SL}_{d+1}$. The linearizations

$$\underline{c} \in \Delta(d+1, n) = \{(c_1, \ldots, c_n) \in \mathbb{Q}^n \mid 0 \leq c_i \leq 1, \sum c_i = d+1\}$$

are parameterized by a hypersimplex. Generators for this ring are given by tableau functions, which appear in many areas of mathematics, particularly representation theory and Schubert calculus [Ful97]. The ideal of relations for $d = 1$ has recently been determined in a landmark series of papers [HMSV09]. The first main result in this paper relates the invariants for $d > 1$ with those of smaller dimension:

**Theorem 1.1.** *Fix integers $n = n_1 + n_2$ and $d = d_1 + d_2$ with $n_i \geq 2, d_i \geq 1$, and consider two generic linear subspaces $\mathbb{P}^{d_i} \subseteq \mathbb{P}^d$. For any $\underline{c} \in \Delta(d+1, n)$ such that*

$$d_1 \leq \sum_{i=1}^{n_1} c_i \leq d_1 + 1 \quad and \quad d_2 \leq \sum_{i=n_1+1}^{n} c_i \leq d_2 + 1,$$

*the restriction of a multidegree $\underline{c}$ invariant to $\mathbb{P}^{d_1} \cup \mathbb{P}^{d_2} \subseteq \mathbb{P}^d$ is the tensor product of multidegree $\underline{c}' := (c_1, \ldots, c_{n_1}, (\sum_{i=n_1+1}^{n} c_i) - d_2)$ and $\underline{c}'' := (c_{n_1+1}, \ldots, c_n, (\sum_{i=1}^{n_1} c_i) - d_1)$ invariants for the $\mathbb{P}^{d_i}$.*

By studying the projective embedding induced by these invariants, we prove:

**Theorem 1.2.** *With notation as above, there is an "attaching" morphism*

$$\gamma : (\mathbb{P}^{d_1})^{n_1+1} /\!\!/_{\underline{c}'} \, \mathrm{SL}_{d_1+1} \times (\mathbb{P}^{d_2})^{n_2+1} /\!\!/_{\underline{c}''} \, \mathrm{SL}_{d_2+1} \longrightarrow (\mathbb{P}^d)^n /\!\!/_{\underline{c}} \, \mathrm{SL}_{d+1} \,.$$

*The projective embedding of the codomain restricts to the Segre product of the projective embeddings of the domain quotients:* $\gamma^* \mathcal{O}_{\underline{c}}(1) \cong \mathcal{O}_{\underline{c}'}(1) \boxtimes \mathcal{O}_{\underline{c}''}(1)$.

These GIT quotients have played an important role in a variety of subjects. For instance, configurations on the line are involved in deep arithmetic constructions [DM86] and configurations in the plane lead naturally to moduli spaces of del Pezzo surfaces [DO88, Dor04]. By duality these quotients parameterize hyperplane arrangements, and studying the modular compactifications one obtains in this context has revealed new insight into the minimal model program [HKT06]. Another ubiquitous compactification of the configuration space of points on the line is the moduli space of stable rational $n$-pointed curves, $\overline{M}_{0,n}$. The interplay between this space and the GIT quotients $(\mathbb{P}^1)^n /\!\!/ \mathrm{SL}_2$ has been studied in [Kap93, AS08, Bol11]

By [Gia11, Theorem 1.1], for any $1 \leq d \leq n-3$ and $\underline{c} \in \Delta(d+1, n)$ there is a morphism $\varphi : \overline{M}_{0,n} \to (\mathbb{P}^d)^n /\!\!/_{\underline{c}} \mathrm{SL}_{d+1}$ sending a configuration of distinct points on $\mathbb{P}^1$ to the corresponding configuration under the $d^{\text{th}}$ Veronese map. The generic point of a boundary divisor $\overline{M}_{0,n_1+1} \times \overline{M}_{0,n_2+1} \subseteq \overline{M}_{0,n}$ gets sent to a configuration supported on a union of two rational normal curves $C_1 \cup C_2 \subseteq \mathbb{P}^{d_1} \cup \mathbb{P}^{d_2}$ with $\deg(C_i) = d_i$ and a node at the point $\mathbb{P}^{d_1} \cap \mathbb{P}^{d_2}$. We define a *GIT bundle* as a line bundle on $\overline{M}_{0,n}$ of the form $\mathcal{G}_{d,\underline{c}} := \varphi^* \mathcal{O}(1)$, and from Theorem 1.2 we deduce:

**Corollary 1.3.** *The restriction of a GIT bundle $\mathcal{G}_{d,\underline{c}}$ to any boundary divisor in $\overline{M}_{0,n}$ is of the form $\mathcal{G}_{d_1,\underline{c}'} \boxtimes \mathcal{G}_{d_2,\underline{c}''}$.*

1.2. **Cyclic covers and Hodge bundles.** One can view a configuration of points on the line as the set of branch points for a ramified cover, and in this way moduli of points leads naturally to moduli of positive genus algebraic curves. More specifically, given integers $n \geq 4, r \geq 2$, and $\underline{c} \in \mathbb{Z}^n$ with $c_i \geq 0$ and $r \mid \sum_{i=1}^n c_i$, there is a map $M_{0,n} \to M_g$ sending a configuration $(\mathbb{P}^1, p_1, \ldots, p_n)$ to the degree $r$ cyclic cover ramified over $\sum_{i=1}^n c_i p_i$. For $p_i = [x_i : 1]$, this is the regular model of the function field extension of $\mathbb{C}(x)$ given by $y^r = (x - x_1)^{c_1} \cdots (x - x_n)^{c_n}$. The genus of the cyclic cover is determined by the Riemann-Hurwitz formula:

$$(1) \qquad g = \frac{1}{2}(2 - 2r + \sum_{i=1}^n (r - \gcd(c_i, r)).$$

This map was studied by Fedorchuk in [Fed11], and he shows that it extends to a morphism $f_{\underline{c}, r} : \overline{M}_{0,n} \to \overline{M}_g$. We prove here that the restriction of this morphism to a boundary divisor is essentially a product of morphisms of the same type:

**Theorem 1.4.** *For $n = n_1 + n_2$ with $n_i \geq 2$, there is a commutative diagram*

(2)

*where $\pi$ is the product of forgetful maps, $\rho$ is an attaching map, $s = g - (g_1 + g_2) + 1$, and $\underline{c}' = (c_1, \ldots, c_{n_1}, \sum_{i=n_1+1}^{n} c_i)$, $\underline{c}'' = (c_{n_1+1}, \ldots, c_n, \sum_{i=1}^{n_1} c_i)$.*

Cyclic covers of the line have found a wide variety of applications—for example, in string theory and mirror symmetry [Roh09], arithmetic geometry (see [LM05] and the references therein), and dynamical systems [EKZ10, McM09]. They have also been used to study special subvarieties in the moduli space of abelian varieties [Moo10, MO11] and curves and divisors in the moduli space of curves [Che10, Fed11]. The restriction of the Hodge bundle $\mathbb{E}_g$ on $\overline{M}_g$ to various cyclic cover loci plays an important role in several of these applications. By studying the Hodge bundles in the above theorem, we deduce:

**Corollary 1.5.** *The restriction of $\det f_{\underline{c},r}^* \mathbb{E}_g$ to any boundary divisor in $\overline{M}_{0,n}$ is of the form $\det f_{\underline{c}',r}^* \mathbb{E}_{g_1} \boxtimes \det f_{\underline{c}'',r}^* \mathbb{E}_{g_2}$.*

We note the structural similarity between these determinant line bundles and the GIT line bundles discussed earlier. This is explained in the following framework.

1.3. **Factorization rules.** The Wess-Zumino-Witten (WZW) model provides an important class of 2-dimensional conformal field theories. It associates to each algebraic curve of a given genus a vector space, known as a space of conformal blocks, fitting together to form an algebraic vector bundle on the moduli stack $\overline{\mathcal{M}}_{g,n}$. These conformal block bundles were widely studied in the 90s due to their connections with mathematical physics, algebraic geometry, and representation theory [TUY89, Fal94, Bea96, Uen97, Wit91, Pau96, LS97, Loo05]. Of particular interest was the Verlinde formula, which computes their rank as a function of the discrete parameters involved. Most proofs of the formula are rooted in the factorization rules, describing for instance the vector space over an irreducible curve in terms of the vector spaces over the components of a reducible curve it degenerates to (cf. §2.4). A key insight was that by abstracting the factorization rules to a numerical function encoding the dimensions under this degeneration process, one can form a ring, called the fusion ring, whose representation theory encodes the information of the Verlinde formula and provides an elegant way to access it (cf. [Bea96]).

In recent years there has been a renewed flurry of activity in the WZW model, in large part due to Fakhruddin's formulae for the Chern classes of conformal block bundles [Fak09]. These formulae provide an intersection-theoretic method for determining the isomorphism class of the associated determinant line bundles. This has revealed much insight into the geometry of these objects, particularly in the case of genus zero [Fak09, AGSS10, AGS10, Swi11, Fed11, Gia11, GG12, GJMS12].

We introduce an abstraction of WZW factorization that is in a sense orthogonal to that of fusion rules: the latter encode the rank of conformal block vector bundles under degeneration while discarding any extraneous information, whereas our framework only applies to rank one bundles but it encodes the isomorphism class of these line bundles. We define a *divisorial factorization system* $\mathcal{L}$ to be a collection of line bundle classes on $\overline{M}_{0,n}$, for $3 \leq n < \infty$, such that the restriction of any $L \in \mathcal{L}$ to any boundary divisor $\overline{M}_{0,m+1} \times \overline{M}_{0,n-m+1} \subseteq \overline{M}_{0,n}$ is of the form $L_1 \boxtimes L_2$ for $L_i \in \mathcal{L}$ (see §2.2). Typically, the members of such a system are indexed by vectors assigning some weight datum to each marked point. We define a *divisorial factorization rule* to be a collection of functions that encode how these weights transform upon restriction to the boundary (Definition 2.4).

**Theorem 1.6.** *If two systems admit the same divisorial factorization rule, then their line bundles are isomorphic for all $n \geq 4$ if they are for $n = 4$.*

The GIT linearization $\underline{c}$ provides a notion of weights for the marked points in the context of GIT bundles $\mathcal{G}_{d,\underline{c}}$ on $\overline{M}_{0,n}$, and by Corollary 1.3 these line bundles form a divisorial factorization system. Similarly, the ramification vector for a cyclic cover yields weights for points on the line and Corollary 1.5 implies that the restrictions of the Hodge class to these weighted cyclic cover loci form another factorization system. These Hodge classes are closely related to line bundles studied by Fedorchuk in [Fed11], defined by first taking a $\mu_r$-eigenbundle decomposition of the restricted Hodge bundle and then considering the determinants of these direct summands. In this paper we prove that these line bundles, which we call *cyclic bundles*, also form a divisorial factorization system, and moreover the following remarkable property holds: the divisorial factorization rule for GIT bundles, cyclic bundles, and a certain class of conformal block (CB) line bundles all coincide. Consequently, we obtain a new proof of the following:

**Theorem 1.7.** *Fix $n \geq 4, r \geq 2$, and $\underline{c} \in \mathbb{Z}^n$ with $0 \leq c_i < r$ and $r \mid \sum_{i=1}^n c_i$. The following line bundles on $\overline{M}_{0,n}$ are isomorphic:*

(i) *the level 1, $\mathfrak{sl}_r$ conformal block bundle with fundamental weights $(\omega_{c_1}, \ldots, \omega_{c_n})$;*

(ii) *the pullback of $\mathcal{O}(1)$ along the map $\overline{M}_{0,n} \to (\mathbb{P}^{\frac{|\underline{c}|}{r}-1})^n /\!\!/_{\underline{c}} \mathrm{SL}_{|\underline{c}|/r}$ sending a configuration of points on the line to the corresponding configuration on a rational normal curve;*

(iii) *the $r^{th}$ tensor power of the determinant of the unit character $\mu_r$-eigenbundle of the Hodge bundle on the locus of ramified degree $r$ cyclic covers with branch points weighted by $\underline{c}$.*

Indeed, standard results apply for the case $\overline{M}_{0,4} \cong \mathbb{P}^1$, so we simply apply Theorem 1.6. The identification of CB bundles with GIT bundles was first proven in [GG12, Theorem 3.2], whereas the identification with cyclic bundles was first proven in [Fed11, Theorem 4.5]. In both cases the proof relies on computing the degree of the restriction to rational curves in the boundary known as F-curves. Once the formulae for these degrees are obtained, the problem reduces to a completely numerical/combinatorial one, namely, showing that the three degree formulae coincide. Beyond this previously observed numerical coincidence of divisor classes, our proof using factorization shows that the isomorphisms in Theorem 1.7 reflect a common functoriality in these three constructions.

**Remark 1.8.** One can extend our notion of divisorial factorization to higher rank vector bundles and to positive genus $\overline{M}_{g,n}$ by mimicking the statement of WZW factorization in that setting. It would be interesting to find geometric constructions of vector bundles that factorize and to compare them with conformal block vector bundles, as in the case of Theorem 1.7 for rank one and genus zero.

**Acknowledgements**: We thank D. Abramovich, R. Cavalieri, and G. Ottaviani for their generous assistance with this work, and we thank M. Fedorchuk, A. Gibney, and C. Sorger for their encouragement and stimulating discussions. The second author was supported by the SNF.

## 2. Factorization of line bundles on $\overline{M}_{0,n}$

In this section we introduce the notion of divisorial factorization, prove that it can be used to encode Picard group classes through a simple induction, and provide examples coming from the WZW model.

2.1. **Preliminaries.** For any $I \subseteq \{1, \ldots, n\}$ with $2 \leq |I| \leq n-2$, there is a boundary divisor inclusion map $\partial_I : \overline{M}_{0,I+1} \times \overline{M}_{0,I^c+1} \hookrightarrow \overline{M}_{0,n}$, where $I^c := \{1, \ldots, n\} \setminus I$ is the complementary index set. This yields a restriction homomorphism

$$\partial_I^* : \mathrm{Pic}(\overline{M}_{0,n}) \to \mathrm{Pic}(\overline{M}_{0,I+1} \times \overline{M}_{0,I^c+1}).$$

Another map of Picard groups that we will use throughout is

$$\mathrm{Pic}(\overline{M}_{0,I+1}) \times \mathrm{Pic}(\overline{M}_{0,I^c+1}) \to \mathrm{Pic}(\overline{M}_{0,I+1} \times \overline{M}_{0,I^c+1})$$

defined by $(L, L') \mapsto L \boxtimes L'$. This latter map is an isomorphism, as one can check with the Künneth formula, though we will not need this fact.

**Convention 2.1.** Unless otherwise specified, we assume for notational convenience that an index set $I$ of cardinality $m$ is of the form $I = \{1, \ldots, m\}$.

2.2. **Main definitions.**

**Definition 2.2.** A *divisorial factorization system* is a subset

$$\mathcal{L} \subseteq \bigcup_{n \geq 3} \mathrm{Pic}(\overline{M}_{0,n})$$

that is closed under boundary restriction: if $L \in \mathcal{L}$ is in $\mathrm{Pic}(\overline{M}_{0,n})$ for some $n \geq 4$, and $I \subseteq \{1, \ldots, n\}$ satisfies $2 \leq |I| \leq n-2$, then $\partial_I^* L = L' \boxtimes L''$ for some $L', L'' \in \mathcal{L}$. Given a set $S$, we say that the system $\mathcal{L}$ is *S-weighted* if there exists

$$\Phi = (\Phi_3, \Phi_4, \ldots) \in \prod_{n \geq 3} \mathrm{Hom}_{\mathrm{Set}}(S^n, \mathrm{Pic}(\overline{M}_{0,n}))$$

such that $\mathcal{L} = \bigcup_{n \geq 3} \mathrm{Im}\, \Phi_n$. We then call $\Phi$ an *S-weighting* for $\mathcal{L}$.

**Remark 2.3.** To define an *S*-weighting on a system, we will frequently describe the maps $\Phi_n : S^n \to \mathrm{Pic}(\overline{M}_{0,n})$ only for certain vectors $\underline{s} \in S^n$; as such, we implicitly take the trivial line bundle class as the image of all other vectors.

The interpretation of this definition is that an *S*-weighting provides a way of specifying a line bundle on $\overline{M}_{0,n}$ by specifying an element of $S$ for each marked point. The restriction of this line bundle to a boundary divisor is then determined by an element of $S$ for each marked point as before, plus elements of $S$ for each of the two attaching points. Note, however, that we do not require the weights determining this restricted line bundle to be unique.

The following definition provides a method for keeping track of the weighting data with respect to boundary restriction.

**Definition 2.4.** Given an *S*-weighted system $(\mathcal{L}, \Phi)$, a *divisorial factorization rule* is a collection of pairs of maps $(\phi_I, \psi_I)$, for each $n \geq 4$ and $I$ as above, where $\phi_I : S^n \to S^{|I|+1}$ and $\psi_I : S^n \to S^{|I^c|+1}$, satisfying

$$\partial_I^* \circ \Phi_n = (\Phi_{|I|+1} \circ \phi_I) \boxtimes (\Phi_{|I^c|+1} \circ \psi_I).$$

2.3. **Inductive structure.**

**Theorem 2.5.** *If two S-weighted systems $(\mathcal{L}, \Phi)$ and $(\mathcal{L}', \Phi')$ admit the same divisorial factorization rule, and $\Phi_4 = \Phi'_4$, then $\Phi_n = \Phi'_n$ for every $n \geq 4$.*

*Proof.* Fix $\underline{s} \in S^n$. In [Kee92] Keel proved that the Chow groups of $\overline{M}_{0,n}$ are generated by the boundary strata and that rational equivalence coincides with numerical equivalence, so two line bundles are isomorphic if they have the same degree on every F-curve $C \subseteq \overline{M}_{0,n}$. Since $C$ is an intersection of boundary divisors, the fact that $\mathcal{L}$ and $\mathcal{L}'$ are closed under boundary restriction implies that $\Phi_n(\underline{s})|_C = \Phi_4(\underline{u})$ for some $\underline{u} \in S^4$, and similarly $\Phi'_n(\underline{s})|_C = \Phi'_4(\underline{v})$, $\underline{v} \in S^4$. The hypothesis that $\Phi$ and $\Phi'$ admit the same factorization rule implies that we can choose $\underline{u} = \underline{v}$, so $\Phi_n(\underline{s})|_C = \Phi'_n(\underline{s})|_C$ by the assumption about $n = 4$. □

2.4. **Examples from WZW.** For any simple complex Lie algebra $\mathfrak{g}$, level $\ell \in \mathbb{N}$, and $n$-tuple $\underline{\lambda}$ of dominant integral weights of level $\ell$, the WZW model produces an algebraic vector bundle on $\overline{\mathcal{M}}_{g,n}$, which we denote by $\mathbb{V}(\mathfrak{g}, \ell, \underline{\lambda})$. The fiber of this bundle over a pointed curve $[C, p_1, \ldots, p_n] \in \overline{\mathcal{M}}_{g,n}$, called a space of *conformal blocks*, is a finite-dimensional complex vector space that can be realized as a space of generalized parabolic theta functions if $C$ is smooth, and which in general is constructed as a space of covariants [LS97, Pau96, TUY89, Uen97, Wit91].

The rank of these bundles is computed by the Verlinde formula, though in many cases it is more convenient to use factorization, which we state here for genus $g = 0$:

**Theorem 2.6.** [Fak09, Proposition 2.4] *For any $\mathfrak{g}, \ell$, and weights $\lambda_1, \ldots, \lambda_n$, there is a natural isomorphism*

$$\partial_I^* \mathbb{V}(\mathfrak{g}, \ell, (\lambda_1, \ldots, \lambda_n)) \cong \bigoplus_{\mu \in P_\ell} \mathbb{V}(\mathfrak{g}, \ell, (\lambda_1, \ldots, \lambda_m, \mu)) \boxtimes \mathbb{V}(\mathfrak{g}, \ell, (\lambda_{m+1}, \ldots, \lambda_n, \mu^*))$$

*where $P_\ell$ is the set of weights of level $\ell$ and $*$ is the natural involution on this set.*

**Example 2.7.** For $\ell = 1$ and $\mathfrak{g}$ simply laced, the conformal block bundles of this type form a divisorial factorization system. Indeed, in this situation $\mathbb{V}(\mathfrak{g}, 1, \underline{\lambda})$ has rank at most one [Fak09, §5], so the direct sum in Theorem 2.6 has at most one nonzero term, thus a CB line bundle pulls back to a product of CB line bundles.

**Example 2.8.** To illustrate weight data and factorization rules, let us fix $\ell = 1$ and $\mathfrak{g} = \mathfrak{sl}_r$. The level 1 dominant integral weights are then the fundamental weights $\omega_i$, $1 \leq i \leq r - 1$, and the zero weight $\omega_0$ corresponding to the trivial representation. We also set $\omega_r := \omega_0$. With this convention, the involution is given by $\omega_i^* = \omega_{r-i}$. The rank of $\mathbb{V}(\mathfrak{sl}_r, 1, (\omega_{c_1}, \ldots, \omega_{c_n}))$ is $\leq 1$, and it is nonzero if and only if the sum of the weights lies in the root lattice: $\sum_{i=1}^n c_i \in r\mathbb{Z}$. From this description it is easy to see that the following defines a $\mathbb{Z}/r\mathbb{Z}$-weighted factorization system and divisorial factorization rule:

$$\Phi_n(c_1, \ldots, c_n) := \mathbb{V}(\mathfrak{sl}_r, 1, (\omega_{c_1}, \ldots, \omega_{c_n})) \text{ for } \sum_{i=1}^n c_i \in r\mathbb{Z},$$

$$\phi_I(c_1, \ldots, c_n) := (c_1, \ldots, c_m, \sum_{i=m+1}^n c_i \mod r),$$

$$\psi_I(c_1, \ldots, c_n) := (c_{m+1}, \ldots, c_n, \sum_{i=1}^m c_i \mod r).$$

## 3. Point configurations

Here we prove Theorems 1.1 and 1.2. The main idea is to embed the GIT quotients in one projective space and then compare the invariant functions that manifest these embeddings. After concluding the proof, we deduce Corollary 1.3.

3.1. **Stability.** Recall that a fractional linearization for the diagonal $\mathrm{SL}_{d+1}$-action on $(\mathbb{P}^d)^n$ is the choice of an ample class $\underline{c} \in \mathbb{Q}_{>0}^n \subseteq \mathbb{Q}^n \cong \mathrm{Pic}((\mathbb{P}^d)^n)_{\mathbb{Q}}$. We think of this as assigning a positive weight to each point in the configuration. It is customary to rescale so that $|\underline{c}| = d + 1$. In this case, a configuration $(p_1, \ldots, p_n)$, $p_i \in \mathbb{P}^d$, is semistable if and only if, for every linear subspace $W \subseteq \mathbb{P}^d$, the inequality $\sum_{p_i \in W} c_i \leq \dim W + 1$ holds; the configuration is stable if and only if this inequality is strict [DH98, Example 3.3.24].

The closure of the space of linearizations with nonempty semistable locus is the hypersimplex $\Delta(d + 1, n)$. Taking the closure allows weights to be zero. If $c_i = 0$ for some $i$, then the GIT quotient is defined by first applying the basepoint-free complete linear system $\mathcal{O}(\underline{c})$ on $(\mathbb{P}^d)^n$, which contracts the $i^{\mathrm{th}}$ copy of $\mathbb{P}^d$: points with zero weight are "forgotten" in the quotient.

Observe that the vectors $\underline{c}'$ and $\underline{c}''$ in Theorem 1.2 are scaled appropriately:

$$|\underline{c}'| = (\sum_{i=1}^{n_1} c_i) + (\sum_{i=n_1+1}^{n} c_i) - d_2 = |\underline{c}| - d_2 = (d+1) - d_2 = d_1 + 1,$$

and similarly $|\underline{c}''| = d_2 + 1$. Moreover, the bounds on $\sum_{i=1}^{n_1} c_i$ and $\sum_{i=n_1+1}^{n} c_i$ ensure that $\underline{c}' \in \Delta(d_1 + 1, n_1 + 1)$ and $\underline{c}'' \in \Delta(d_2 + 1, n_2 + 1)$.

3.2. **Attaching map.** Let $U \subseteq (\mathbb{P}^{d_1})^{n_1+1} \times (\mathbb{P}^{d_2})^{n_2+1}$ be the closed subset defined by the coordinate matrices (3). Consider the map $\widetilde{\gamma} : U \to (\mathbb{P}^d)^n$ defined as follows:

$$
(3) \qquad (
\begin{bmatrix}
x_{01} & \cdots & x_{0n_1} & 0 \\
\vdots & & \vdots & \vdots \\
x_{d_1-1,1} & \cdots & x_{d_1-1,n_1} & 0 \\
x_{d_1 1} & \cdots & x_{d_1 n_1} & 1
\end{bmatrix},
\begin{bmatrix}
x_{0,n_1+1} & \cdots & x_{0n} & 1 \\
x_{1,n_1+1} & \cdots & x_{1n} & 0 \\
\vdots & & \vdots & \vdots \\
x_{d_2,n_1+1} & \cdots & x_{d_2 n} & 0
\end{bmatrix})
$$

$$
(4) \qquad \mapsto
\begin{bmatrix}
x_{01} & \cdots & x_{0n_1} & 0 & \cdots & 0 \\
\vdots & & \vdots & \vdots & & \vdots \\
x_{d_1-1,1} & \cdots & x_{d_1-1,n_1} & 0 & \cdots & 0 \\
x_{d_1 1} & \cdots & x_{d_1 n_1} & x_{0,n_1+1} & \cdots & x_{0n} \\
0 & \cdots & 0 & x_{1,n_1+1} & \cdots & x_{1n} \\
\vdots & & \vdots & \vdots & & \vdots \\
0 & \cdots & 0 & x_{d_2,n_1+1} & \cdots & x_{d_2 n}
\end{bmatrix}.
$$

Geometrically, $\widetilde{\gamma}$ uses the fixed $(n_i + 1)^{\mathrm{th}}$ point in $\mathbb{P}^{d_i}$ to attach these two subspaces and embed them in $\mathbb{P}^d$, where the attaching point then has coordinates

$$q := [0 : \cdots : 0 : 1 : 0 : \cdots : 0] = \mathbb{P}^{d_1} \cap \mathbb{P}^{d_2} \subseteq \mathbb{P}^d.$$

**Lemma 3.1.** *If $U_{ss} := U \cap ((\mathbb{P}^{d_1})_{ss}^{n_1+1} \times (\mathbb{P}^{d_2})_{ss}^{n_2+1})$ denotes the semistable locus with respect to both $\mathrm{SL}_{d_i+1}$ actions, then $\widetilde{\gamma}(U_{ss}) \subseteq (\mathbb{P}^d)_{ss}^n$.*

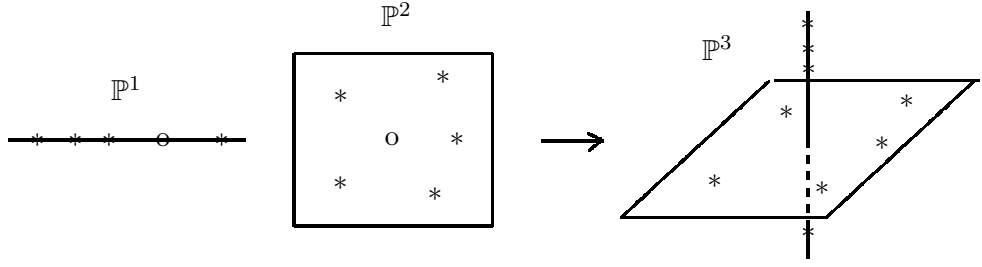

FIGURE 1. The map $\widetilde{\gamma} : U \to (\mathbb{P}^3)^9$, where $U \subseteq (\mathbb{P}^1)^{4+1} \times (\mathbb{P}^2)^{5+1}$ denotes the set of configurations such that the attaching points have coordinates $[0:1]$ and $[1:0:0]$, respectively.

*Proof.* Let $W \subseteq \mathbb{P}^d$ be a linear subspace, and first consider the case $q \notin W$. Then

$$\dim W = \dim(W \cap \mathbb{P}^{d_1}) + \dim(W \cap \mathbb{P}^{d_2}) + 1,$$

where by convention $\dim \varnothing := -1$. Therefore,

$$\sum_{p_i \in W} c_i = \sum_{p_i \in W \cap \mathbb{P}^{d_1}} c_i + \sum_{p_i \in W \cap \mathbb{P}^{d_2}} c_i \leq (\dim(W \cap \mathbb{P}^{d_1}) + 1) + (\dim(W \cap \mathbb{P}^{d_2}) + 1) = \dim W + 1,$$

as required.

Next, suppose that $q \in W$. Then

$$\sum_{p_i \in W \cap \mathbb{P}^{d_1}} c_i \leq \dim(W \cap \mathbb{P}^{d_1}) + 1 - ((\sum_{i=n_1+1}^{n} c_i) - d_2)$$

because a semistable configuration in $(\mathbb{P}^{d_1})^{n_1+1}$ has weight $\leq \dim(W \cap \mathbb{P}^{d_1}) + 1$ in the linear subspace $W \cap \mathbb{P}^{d_1}$, and applying $\widetilde{\gamma}$ has the effect of dropping the extra point of weight $(\sum_{i=n_1+1}^{n} c_i) - d_2$. The analogous inequality holds for $W \cap \mathbb{P}^{d_2}$, so

$$\sum_{p_i \in W} c_i \leq \sum_{p_i \in W \cap \mathbb{P}^{d_1}} c_i + \sum_{p_i \in W \cap \mathbb{P}^{d_2}} c_i \leq \dim(W \cap \mathbb{P}^{d_1}) + \dim(W \cap \mathbb{P}^{d_2}) + 2 + (d_1 + d_2) - \sum_{i=1}^{n} c_i.$$

But $d_1 + d_2 = d$, $\sum_{i=1}^{n} c_i = d + 1$, and $\dim(W \cap \mathbb{P}^{d_1}) + \dim(W \cap \mathbb{P}^{d_2}) = \dim W$, so the right side of the preceding inequality reduces to $\dim W + 1$, as required.    □

3.3. **Invariants and embeddings.** We may assume that the linearization vectors in Theorem 1.2 are integral, since the hypothesis and conclusion are invariant under scaling. Specifically, choose $k \in \mathbb{Z}$ such that $k\underline{c} \in \mathbb{Z}^n$, and write $\underline{k} := k\underline{c}$, so that $k_i = kc_i$. The condition $\underline{c} \in \Delta(d+1,n)$ implies that $|\underline{k}| = k(d+1)$, which in turn implies that there exist invariant functions on $(\mathbb{P}^d)^n$ of multi-degree $\underline{k}$. Indeed, it is well-known [Dol03, §11.2] that the invariant functions

$$T := \mathrm{H}^0((\mathbb{P}^d)^n, \mathcal{O}(\underline{k}))^{\mathrm{SL}_{d+1}} \cong (S^{k_1}\mathbb{C}^{d+1} \otimes \cdots \otimes S^{k_n}\mathbb{C}^{d+1})^{\mathrm{SL}_{d+1}}$$

generate the ring of $\mathrm{SL}_{d+1}$-invariants and this linear system induces the composition

$$(\mathbb{P}^d)^n_{ss} \to (\mathbb{P}^d)^n /\!\!/_{\underline{k}} \mathrm{SL}_{d+1} \hookrightarrow \mathbb{P}(T^*).$$

If we choose a $(d+1) \times n$ matrix of homogeneous coordinates for $(\mathbb{P}^d)^n$, then a basis for $T$ is given by the multi-degree $\underline{k}$ products of minors of this matrix. More precisely, a basis is the set of semistandard tableau functions of size $(d +$ $1) \times k$ with entries in $\{1, \ldots, n\}$ such that $i$ appears $k_i$ times. Each column of a tableau represents the minor of the coordinate matrix determined by the numbers appearing in the tableau column; juxtaposition of columns indicates the product of the corresponding minors. The semistandard condition means that the entries decrease weakly across the rows and strongly down the columns. See [Dol03, §5.6], [Muk03, §8.1(C)], and [Ful97] for more background.

Define $\underline{k}' := k\underline{c}'$, $T_1 := \mathrm{H}^0((\mathbb{P}^{d_1})^{n_1+1}, \mathcal{O}(\underline{k}'))^{\mathrm{SL}_{d_1+1}}$, and similarly for $\underline{k}''$, $T_2$. Then the linear system $T_1 \boxtimes T_2$ induces the composition

$$(\mathbb{P}^{d_1})_{ss}^{n_1+1} \times (\mathbb{P}^{d_2})_{ss}^{n_2+1} \to (\mathbb{P}^{d_1})^{n_1+1} /\!\!/_{\underline{k}'} \mathrm{SL}_{d_1+1} \times (\mathbb{P}^{d_2})^{n_2+1} /\!\!/_{\underline{k}''} \mathrm{SL}_{d_2+1}$$
$$\hookrightarrow \mathbb{P}(T_1^*) \times \mathbb{P}(T_2^*) \hookrightarrow \mathbb{P}(T_1^* \otimes T_2^*).$$

Combining this with the attaching map defined above, we obtain the diagram

$$(5) \quad \begin{array}{ccccc} U_{ss} & \longrightarrow & (\mathbb{P}^{d_1})^{n_1+1} /\!\!/_{\underline{k}'} \mathrm{SL}_{d_1+1} \times (\mathbb{P}^{d_2})^{n_2+1} /\!\!/_{\underline{k}''} \mathrm{SL}_{d_2+1} & \hookrightarrow & \mathbb{P}(T_1^* \otimes T_2^*) \\ \widetilde{\gamma} \downarrow & & & & \\ (\mathbb{P}^d)_{ss}^n & \longrightarrow & (\mathbb{P}^d)^n /\!\!/_{\underline{k}} \mathrm{SL}_{d+1} & \hookrightarrow & \mathbb{P}(T^*) \end{array}$$

where surjectivity of the top-left arrow follows from the fact that every orbit for the $\mathrm{SL}_{d_1+1} \times \mathrm{SL}_{d_2+1}$-action has a representative in $U$.

Let $T_1 \otimes T_2|_U := \{f|_U : f \in T_1 \otimes T_2\}$ and similarly for $T|_{\widetilde{\gamma}(U)}$. Note that $T_1 \otimes T_2|_U \cong T_1 \otimes T_2$, since all these functions are invariant and their orbits are represented in $U$. The following is a refinement of, and immediately implies, Theorem 1.1:

**Proposition 3.2.** *There is a surjection $T|_{\widetilde{\gamma}(U)} \twoheadrightarrow T_1 \otimes T_2|_U$. The composition*

$$T \twoheadrightarrow T|_{\widetilde{\gamma}(U)} \twoheadrightarrow T_1 \otimes T_2|_U \xrightarrow{\sim} T_1 \otimes T_2$$

*induces a linear embedding $\mathbb{P}(T_1^* \otimes T_2^*) \hookrightarrow \mathbb{P}(T^*)$ making diagram (5) commute.*

*Proof.* Denote by $A_1$ and $A_2$ the matrices in (3) parameterizing $U$, and number their columns $(1, \ldots, n_1 + 1)$ and $(1, \ldots, n_2 + 1)$, respectively. Denote by $B$ the matrix in (4) parameterizing $\widetilde{\gamma}(U)$, with columns $(1, \ldots, n)$. We call the first $n_1$ columns of $B$ the *first block*, and the next $n_2$ columns the *second block*.

Given a $(d+1) \times 1$ tableau function on $(\mathbb{P}^d)^n$, its restriction to $\widetilde{\gamma}(U)$ is zero unless its entries specify $d_1 + 1$ columns from the first block and $d_2$ from the second block, or $d_1$ from the first and $d_2 + 1$ from the second. In this case, the restriction is the product of a $(d_1 + 1) \times 1$ tableau on $A_1$ and a $(d_2 + 1) \times 1$ tableau on $A_2$. The entry $n_2 + 1$ appears exactly once in the $A_2$ tableau if the original tableau specifies $d_1 + 1$ columns from the first block and $d_2$ from the second, and vice versa for $n_1 + 1$. For example, if

$$f = \begin{array}{|c|} \hline 1 \\ \hline \vdots \\ \hline d_1 + 1 \\ \hline n_1 + 1 \\ \hline \vdots \\ \hline n_1 + d_2 \\ \hline \end{array} \text{, then } f|_{\widetilde{\gamma}(U)} = \begin{array}{|c|} \hline 1 \\ \hline \vdots \\ \hline d_1 + 1 \\ \hline \end{array} \otimes \begin{array}{|c|} \hline 1 \\ \hline \vdots \\ \hline d_2 \\ \hline n_2 + 1 \\ \hline \end{array}$$

As mentioned above, a basis for $T$ is given by certain $(d+1) \times k$ tableau functions. Since these are a $k$-fold product of $(d+1) \times 1$ tableau functions, their restriction to $\widetilde{\gamma}(U)$ is zero unless the entries of each column are distributed across the columns of $B$ as in the preceding paragraph. Therefore, there is a linear map

$$\mu : T \to \bigoplus_{(\underline{m}_1, \underline{m}_2) \in \mathbb{N}^{n+2}} \mathrm{H}^0(U, \mathcal{O}(\underline{m}_1) \boxtimes \mathcal{O}(\underline{m}_2)|_U)$$

defined by first restricting each basis element to $\widetilde{\gamma}(U)$ and then expressing as a tensor product of a $(d_1 + 1) \times k$ tableau on $A_1$ and a $(d_2 + 1) \times k$ tableau on $A_2$. We claim that the image of $\mu$ is $T_1 \otimes T_2|_U$.

Since $|\underline{k}'| = k|\underline{c}'| = k(d_1 + 1)$ and $|\underline{k}''| = k(d_2 + 1)$, the image under $\mu$ of a basis tableau of $T$ (with nonzero restriction) has the same dimensions as the basis elements of $T_1 \otimes T_2|_U$. Moreover, the fact that the original tableau is semistandard immediately implies that the tableau functions in the restriction are semistandard with entries in $\{1, \ldots, n_1 + 1\}$ and $\{1, \ldots, n_2 + 1\}$.

We will show that:

(1) the index $i$ occurs exactly $kc_i$ times on the $A_1$ factor, for $i \le n_1$, and $kc_{n_1+i}$ times on the $A_2$ factor, for $i \le n_2$;
(2) $n_1 + 1$ occurs $k(\sum_{i=n_1+1}^{n} c_i) - kd_2$ times on the $A_1$ factor;
(3) $n_2 + 1$ occurs $k(\sum_{i=1}^{n_1} c_i) - kd_1$ times on the $A_2$ factor.

By the definition of $\underline{k}'$ and $\underline{k}''$, this is the form of the basis elements of $T_1 \otimes T_2|_U$, so this will prove that $\mu(T) \subseteq T_1 \otimes T_2|_U$.

If $i \le n_1$, then $i$ occurs $kc_i$ times in the original tableau, and by our explicit description above we see that $i$ occurs $kc_i$ times in the restriction. Similarly, if $i \ge n_1 + 1$ then it occurs $kc_{i-n_1}$ times in the original and hence $kc_i$ times in the restriction. This verifies part (1) of the claim.

By the nonzero restriction to $U$ assumption, each column of the original $(d+1) \times k$ tableau has $d_1 + 1$ entries from the first block and $d_2$ from the second, or $d_1$ from the first and $d_2 + 1$ from the second. We saw that the index $n_2 + 1$ occurs once in the restriction in the first case and $n_1 + 1$ occurs once in the second case. If we denote by $\alpha$ the number of columns of the original tableau of the first type, and by $\beta$ the number of the second, then the multiplicity of $n_2 + 1$ is $\alpha$ and the multiplicity of $n_1 + 1$ is $\beta$. The total number of entries in the tableau corresponding to columns in the first block is by assumption $\sum_{i=1}^{n_1} kc_i$. On the other hand, this number is also $\alpha(d_1 + 1) + \beta d_1$, since $\alpha$ columns of the tableau have $d_1 + 1$ such entries and $\beta$ columns have $d_1$ such entries. Clearly $\alpha + \beta = k$, the total number of columns, so

$$\sum_{i=1}^{n_1} kc_i = \alpha(d_1 + 1) + \beta d_1 = \alpha + kd_1,$$

which implies that $n_2 + 1$ indeed has multiplicity $\alpha = k(\sum_{i=1}^{n_1} c_i) - kd_1$. Analogously, one sees that $n_1 + 1$ has multiplicity $\beta = k(\sum_{i=n_1+1}^{n} c_i) - kd_2$, verifying parts (2) and (3) of the claim.

Moreover, it follows from this discussion that the containment $\mu(T) \subseteq T_1 \otimes T_2|_U$ is in fact an equality, since any pair of semistandard tableau functions with entries as in parts (1), (2) and (3) of the claim is the restriction of a basis element of $T$.

Finally, the fact that the map $\mathbb{P}(T_1^* \otimes T_2^*) \hookrightarrow \mathbb{P}(T^*)$ induced by $\mu$ renders diagram (5) commutative is clear from construction. □

3.4. **Concluding the proof of Theorem 1.2.** Proposition 3.2 yields an embedding

$$(\mathbb{P}^{d_1})^{n_1+1} /\!\!/_{\underline{k}'} \mathrm{SL}_{d_1+1} \times (\mathbb{P}^{d_2})^{n_2+1} /\!\!/_{\underline{k}''} \mathrm{SL}_{d_2+1} \hookrightarrow \mathbb{P}(T^*),$$

and it follows from commutativity of diagram (5), appended by the map in Proposition 3.2, and surjectivity of the horizontal map from $U_{ss}$ there that the image of this embedding is contained in the image of the embedding $(\mathbb{P}^d)^n /\!\!/_{\underline{k}} \mathrm{SL}_{d+1} \hookrightarrow \mathbb{P}(T^*)$. Therefore, we obtained the desired morphism $\gamma$ in the statement of Theorem 1.2:

$$
\begin{array}{ccc}
(\mathbb{P}^{d_1})^{n_1+1} /\!\!/_{\underline{k}'} \mathrm{SL}_{d_1+1} \times (\mathbb{P}^{d_2})^{n_2+1} /\!\!/_{\underline{k}''} \mathrm{SL}_{d_2+1} & \hookrightarrow & \mathbb{P}(T_1^* \otimes T_2^*) \\
\downarrow{\scriptstyle \gamma} & & \uparrow \\
(\mathbb{P}^d)^n /\!\!/_{\underline{k}} \mathrm{SL}_{d+1} & \hookrightarrow & \mathbb{P}(T^*)
\end{array}
$$

Moreover, the statement about polarizations follows immediately since the inclusion $\mathbb{P}(T_1^* \otimes T_2^*) \hookrightarrow \mathbb{P}(T^*)$ is linear and all the GIT quotients involved inherit their hyperplane class from their embeddings in these projective spaces.    □

3.5. **Factorization.** We are now in a position to deduce the following, which is an explicit version of Corollary 1.3 stated in the introduction:

**Corollary 3.3.** *Let* $\mathcal{G}_{d,\underline{c}} := \varphi_{d,\underline{c}}^* \mathcal{O}(1)$, *where* $\varphi_{d,\underline{c}} : \overline{M}_{0,n} \to (\mathbb{P}^d)^n /\!\!/_{\underline{c}} \mathrm{SL}_{d+1}$ *and* $\mathcal{O}(1)$ *is the GIT polarization. Then, with notation as in Theorem 1.2 and §2.1,*

$$\partial_I^* \mathcal{G}_{d,\underline{c}} \cong \mathcal{G}_{d_1,\underline{c}'} \boxtimes \mathcal{G}_{d_2,\underline{c}''}.$$

*Proof.* By Theorem 1.2, it is enough to show that the following is commutative:

$$
\begin{array}{ccc}
\overline{M}_{0,n_1+1} \times \overline{M}_{0,n_2+1} & \hookrightarrow & \overline{M}_{0,n} \\
\downarrow{\scriptstyle (\varphi_{d_1,\underline{c}'}, \varphi_{d_2,\underline{c}''})} & & \downarrow{\scriptstyle \varphi_{d,\underline{c}}} \\
(\mathbb{P}^{d_1})^{n_1+1} /\!\!/_{\underline{c}'} \mathrm{SL}_{d_1+1} \times (\mathbb{P}^{d_2})^{n_2+1} /\!\!/_{\underline{c}''} \mathrm{SL}_{d_2+1} & \longrightarrow & (\mathbb{P}^d)^n /\!\!/_{\underline{c}} \mathrm{SL}_{d+1}
\end{array}
$$

To check that this diagram commutes, it is enough to consider the restriction to $M_{0,n_1+1} \times M_{0,n_2+1}$, since all the varieties involved are separated. On this open locus, traversing the diagram in either direction corresponds to sending a nodal curve $C = C_1 \cup C_2$ to a union of two rational normal curves in $\mathbb{P}^{d_1} \cup \mathbb{P}^{d_2} \subseteq \mathbb{P}^d$ with a node at the attaching point $\mathbb{P}^{d_1} \cap \mathbb{P}^{d_2}$. It follows from the explicit description of $\gamma$ and of the morphisms in [Gia11, §4.1] that the resulting configurations are projectively equivalent.    □

## 4. Cyclic covers

In this section we prove Theorem 1.4 and Corollary 1.5 and discuss the eigenbundle determinants studied by Fedorchuk.

4.1. **Attaching maps and Hodge bundles.** Let us denote the Hodge bundle over $\overline{M}_g$ by $\mathbb{E}_g$. For any integers $g_1, g_2 \geq 0$ and $s \geq 1$, there is a natural attaching map $\rho : \overline{M}_{g_1,s} \times \overline{M}_{g_2,s} \to \overline{M}_g$, where $g = g_1 + g_2 + s - 1$, defined by glueing two curves of genus $g_1, g_2$ along their $s$ marked points. See Figure 2 for an example. We will also consider the product of forgetful maps $\pi : \overline{M}_{g_1,s} \times \overline{M}_{g_2,s} \to \overline{M}_{g_1} \times \overline{M}_{g_2}$.

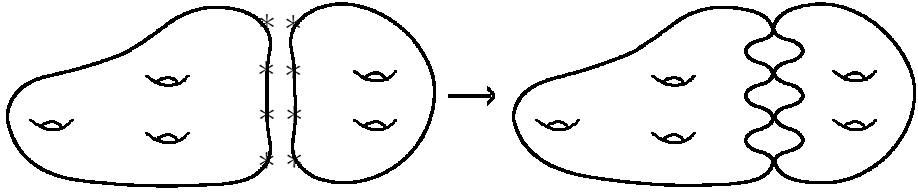

FIGURE 2. The attaching map $\rho : \overline{M}_{3,4} \times \overline{M}_{2,4} \to \overline{M}_8$.

**Lemma 4.1.** *There is an isomorphism of vector bundles*

$$\rho^* \mathbb{E}_g \cong \pi^* ((\mathbb{E}_{g_1} \boxtimes \mathcal{O}_{\overline{M}_{g_2}}) \oplus (\mathcal{O}_{\overline{M}_{g_1}} \boxtimes \mathbb{E}_{g_2})) \oplus \mathcal{O}_{\overline{M}_{g_1,s} \times \overline{M}_{g_2,s}}^{\oplus s-1}.$$

*Proof.* This is standard and has appeared in the literature in various forms. See, e.g., [FP00, Proof of Proposition 2], [CY10, Equation (23)], and [Mum83]. □

4.2. **Admissible covers.** The cyclic cover morphisms $f_{\underline{c},r}$ mentioned in §1.2 are most easily defined on the interior $M_{0,n}$, and Fedorchuk in [Fed11, §2] describes two ways to see the extension to the stable compactification $\overline{M}_{0,n}$. One method is via orbicurves and essentially shows that the extension is obtained by taking a degree $r$ cyclic cover of the universal curve not of $\overline{M}_{0,n}$ itself but of a DM-stack lying over $\overline{M}_{0,n}$, namely Fedorchuk's moduli space of $r$-divisible curves. These $r$-divisible curves are orbicurves whose coarse moduli space is a marked stable rational curve. The second method to see the extension is via the Harris-Mumford theory of admissible covers: to compute the limit of a family of cyclic covers, take the limit in the compactified Hurwitz scheme of admissible covers [HM82]. We will use the latter perspective to prove to Theorem 1.4. To begin, let us see what happens to smooth cyclic covers as they degenerate in the simplest possible way:

**Lemma 4.2.** *Let $C_0 \to D_0 = \mathbb{P}^1 \cup_{\{q\}} \mathbb{P}^1$ be the admissible cover obtained as the limit of cyclic covers $C_t \to \mathbb{P}^1$ ramified over $\sum_{i=1}^n c_i p_i$ as $(\mathbb{P}^1, p_1, \dots, p_n) \in M_{0,n}$ tends to the generic point $(\mathbb{P}^1 \cup_{\{q\}} \mathbb{P}^1, p_1, \dots, p_n)$ of a boundary divisor $\overline{M}_{0,n_1+1} \times \overline{M}_{0,n_2+1} \subseteq \overline{M}_{0,n}$. Let $C' \sqcup C'' \to C_0$ and $D' \sqcup D'' \to D_0$ denote the normalizations, so $D' = \mathbb{P}^1$ with markings $p_1, \dots, p_{n_1}, q$ and $D'' = \mathbb{P}^1$ with markings $p_{n_1+1}, \dots, p_n, q$. Then*

(1) *$C' \to D'$ is a cyclic cover ramified over $\sum_{i=1}^{n_1} c_i p_i + (\sum_{i=n_1+1}^n c_i) q$;*

(2) *$C'' \to D''$ is a cyclic cover ramified over $\sum_{i=n_1+1}^n c_i p_i + (\sum_{i=1}^{n_1} c_i) q$;*

(3) *The fiber of $C_0 \to D_0$ over the node $q$ consists of $s = \gcd(\sum_{i=n_1+1}^n c_i, r) = \gcd(\sum_{i=1}^{n_1} c_i, r)$ distinct points;*

(4) *The genus $g$ of $C_t$ satisfies $g = g_1 + g_2 + s - 1$, where $g_1$ is the genus of $C'$ and $g_2$ is the genus of $C''$.*

*Proof.* By the definition of admissible covers, the restricted maps $C' \to D'$ and $C'' \to D''$ are degree $r$ cyclic covers of $\mathbb{P}^1$. Moreover, they are ramified over a weighted sum of the marked points and the attaching point $q$, with weights on the marked points given by the original weight vector $\underline{c}$. The weight on the attaching point is then uniquely determined, up to equivalence modulo $r$, by the condition that the sum of the weights for each of the two cyclic covers is divisible by $r$. More precisely, since $\sum_{i=1}^n c_i \in r\mathbb{Z}$, it follows immediately that the weights of the attaching points must be as in the statement of the lemma. This proves the first two claims.

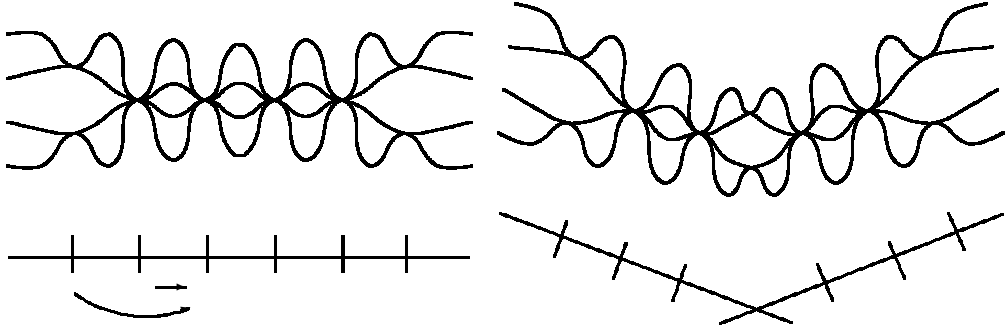

FIGURE 3. The limit of a degree $r = 4$ cover with weights $\underline{c} = (2, 1, 3, 3, 1, 2)$ is obtained by gluing covers with weights $\underline{c}' = (2, 1, 3, 2)$ and $\underline{c}'' = (2, 3, 1, 2)$ at two points lying over the branch node.

The genus of $C_t$ is independent of $t$ and is given by the Riemann-Hurwitz formula:

$$g = 1 - r + \frac{1}{2} \sum_{i=1}^{n} (r - \gcd(c_i, r)).$$

On the other hand, applying Riemann-Hurwitz to $C'$ and $C''$ yields

$$g_1 = 1 - r + \frac{1}{2} \sum_{i=1}^{n_1} (r - \gcd(c_i, r)) + \frac{1}{2}(r - \gcd(\sum_{i=n_1+1}^{n} c_i, r)),$$

$$g_2 = 1 - r + \frac{1}{2} \sum_{i=n_1+1}^{n} (r - \gcd(c_i, r)) + \frac{1}{2}(r - \gcd(\sum_{i=1}^{m} c_i, r)).$$

Since $\gcd(\sum_{i=n_1+1}^{n} c_i, r) = \gcd(\sum_{i=1}^{n_1} c_i, r)$, and we call this common number $s$, which is the number of nodes in the fiber over the node $q$, we see that $g_1 + g_2 = g - s + 1$. This concludes the third and fourth claim. $\qquad\square$

Let $\Delta \subseteq \overline{M}_g$ denote the image of the attaching map $\rho : \overline{M}_{g_1, s} \times \overline{M}_{g_2, s} \to \overline{M}_g$ from §4.1. The normalization of $\Delta$ is $(\overline{M}_{g_1, s} \times \overline{M}_{g_2, s})/S_s$. It follows from Lemma 4.2 that the restriction of the cyclic cover morphism $f_{\underline{c}, r} : \overline{M}_{0, n} \to \overline{M}_g$ to the boundary divisor $\overline{M}_{0, n_1+} \times \overline{M}_{n_2+1}$ has image contained in $\Delta$. Since $\overline{M}_{0, n_1+} \times \overline{M}_{n_2+1}$ is normal, this restriction factors through $(\overline{M}_{g_1, s} \times \overline{M}_{g_2, s})/S_s$. We thus have a diagram as in Theorem 1.4 which is commutative except possibly at the triangle on the left:

However, to check this remaining commutativity we can restrict to the dense open subset $M_{0,n_1+1} \times M_{0,n_2+1} \subseteq \overline{M}_{0,n_1+1} \times \overline{M}_{0,n_2+1}$, where commutativity then follows from Lemma 4.2. This completes the proof of Theorem 1.4.

### 4.3. Determinants and eigenbundles.

Following [Fed11, §4], we define the cyclic bundle $\Sigma_{r,\underline{c}}$ by splitting the Hodge bundle on the cyclic cover locus into $\mu_r$-eigenbundles,

$$f_{\underline{c},r}^* \mathbb{E}_g \cong \bigoplus_{i=1}^{r-1} \mathbb{E}_g^{(i)},$$

and then taking the determinant of the first eigenbundle: $\Sigma_{r,\underline{c}} := \det \mathbb{E}_g^{(1)}$.

**Corollary 4.3.** *With notation as above, we have $\partial_I^* \Sigma_{r,\underline{c}} \cong \Sigma_{r,\underline{c}'} \boxtimes \Sigma_{r,\underline{c}''}$.*

*Proof.* By Theorem 1.4 and Lemma 4.1, the restriction of $f_{\underline{c},r}^* \mathbb{E}_g$ to the boundary divisor $\overline{M}_{0,n_1+1} \times \overline{M}_{0,n_2+1} \subseteq \overline{M}_{0,n}$ coincides with the pullback of $\mathbb{E}_{g_1} \oplus \mathbb{E}_{g_2} \oplus \mathcal{O}^{\oplus s-1}$ along the product morphism $(f_{\underline{c}',r}, f_{\underline{c}'',r}) : \overline{M}_{0,n_1+1} \times \overline{M}_{0,n_2+1} \to \overline{M}_{g_1} \times \overline{M}_{g_2}$. The $\mu_r$-action is compatible with this decomposition, and determinants commute with base change, so

$$\partial_I^* \Sigma_{r,\underline{c}} \cong \det \mathbb{E}_{g_1}^{(1)} \otimes \det \mathbb{E}_{g_2}^{(1)} \otimes \det((\mathcal{O}^{\oplus s-1})^{(1)}).$$

But over a complete variety every direct summand of a trivial bundle is trivial, since endomorphisms are simply constant matrices, so the last factor on the right-hand side vanishes and we obtain the desired result.  □

## 5. Conclusion

Here we use the results from the previous sections to derive Theorem 1.7, restated below for the convenience of the reader:

**Theorem 5.1.** *Fix $n \geq 4, r \geq 2$, and $\underline{c} \in \mathbb{Z}^n$ with $0 \leq c_i < r$ and $r \mid \sum_{i=1}^n c_i$. The following line bundles on $\overline{M}_{0,n}$ are isomorphic:*

(i) *the CB bundle $\mathbb{V}(\mathfrak{sl}_r, 1, (\omega_{c_1}, \ldots, \omega_{c_n}))$;*

(ii) *the GIT bundle $\mathcal{G}_{\frac{|\underline{c}|}{r}-1,\underline{c}}$;*

(iii) *the $r^{th}$ tensor power of the cyclic bundle $\Sigma_{r,\underline{c}}$.*

This is done by showing that all three classes of line bundles form divisorial factorization systems admitting the same factorization rule. Once this is demonstrated, all that remains is to analyze the case $n = 4$ and apply Theorem 2.5.

**Proposition 5.2.** *For each $n, r$ and $\underline{c}$ as above, if $\Phi_n^{CB}(\underline{c}) := \mathbb{V}(\mathfrak{sl}_r, 1, (\omega_{c_1}, \ldots, \omega_{c_n}))$, $\Phi_n^{\mathcal{G}}(\underline{c}) := \mathcal{G}_{\frac{|\underline{c}|}{r}-1,\underline{c}}$, and $\Phi_n^{\Sigma}(\underline{c}) := \Sigma_{r,\underline{c}}^{\otimes r}$, then these three assignments form divisorial factorization systems that all admit the following factorization rule:*

$$\phi_I(c_1, \ldots, c_n) := (c_1, \ldots, c_m, \sum_{i=m+1}^n c_i \mod r)$$

$$\psi_I(c_1, \ldots, c_n) := (c_{m+1}, \ldots, c_n, \sum_{i=1}^m c_i \mod r),$$

*where the mod $r$ representative is taken in $\{1, \ldots, r\}$ for $\phi_I$ and in $\{0, \ldots, r-1\}$ for $\psi_I$. By convention, we set $\mathcal{G}_{d,\underline{c}} = \mathcal{O}$ if $d \notin \{1, \ldots, n-3\}$.*

**Remark 5.3.** By definition, the representative of the remainder is irrelevant for the CB and cyclic bundles, but for GIT we must break the symmetry, as we will see in the proof below.

*Proof.* That this holds for $\Phi^{CB}$ and $\Phi^\Sigma$ was established in Example 2.8 and Corollary 4.3, respectively, so all that remains is establishing it for the GIT system $\Phi^{\mathcal{G}}$. For this, we must show that

$$\partial_I^* \mathcal{G}_{\frac{|\underline{c}|}{r}-1,\underline{c}} \cong \mathcal{G}_{\frac{|\underline{c}'|}{r}-1,\underline{c}'} \boxtimes \mathcal{G}_{\frac{|\underline{c}''|}{r}-1,\underline{c}''},$$

where $\underline{c}' := (c_1, \ldots, c_m, \sum_{i=m+1}^n c_i \mod r)$ and $\underline{c}'' := (c_{m+1}, \ldots, c_n, \sum_{i=1}^m c_i \mod r)$. By rescaling the GIT polarizations, it clearly suffices to prove that

$$\partial_I^* \mathcal{G}_{|\underline{b}|-1,\underline{b}} \cong \mathcal{G}_{|\underline{b}'|-1,\underline{b}'} \boxtimes \mathcal{G}_{|\underline{b}''|-1,\underline{b}''},$$

where $\underline{b} := \frac{1}{r}\underline{c}$, $\underline{b}' := \frac{1}{r}\underline{c}'$, and $\underline{b}'' := \frac{1}{r}\underline{c}''$. If we define $d := |\underline{b}| - 1$, $d_1 := |\underline{b}'| - 1$, and $d_2 := |\underline{b}''| - 1$, then observe that these $d, d_1, d_2$ are all integers. We can assume that they are all positive, for if not then the corresponding GIT bundle is defined to be the trivial line bundle; since any restriction of a trivial bundle is trivial, this is compatible with the above factorization rules. Then $\underline{b} \in \Delta(d+1, n)$, $\underline{b}' \in \Delta(d_1+1, m+1)$, and $\underline{b}'' \in \Delta(d_2+1, n-m+1)$. Our convention on representatives of the remainder modulo $r$ implies that $d_1 + d_2 = d$. The result them follows immediately from Corollary 3.3. □

**5.1. The case $n = 4$.** Since $\overline{M}_{0,4} \cong \mathbb{P}^1$ and $\mathrm{Pic}(\overline{M}_{0,4}) \cong \mathbb{Z}$, the isomorphism class of any line bundle here is simply its degree. To complete the proof of Theorem 1.7, by applying Theorem 2.5, we must show that the three types of line bundles occurring have the same degree for $n = 4$. Although this is by no means trivial, the computation of these degrees has already appeared in the literature:

**Theorem 5.4.** *With notation as in Proposition 5.2, and assuming for notational convenience that $\underline{c} \in \mathbb{Z}^4$ satisfies $c_1 \leq \cdots \leq c_4 \leq r$, we have:*

- $\deg \Phi_4^{CB}(\underline{c}) = \begin{cases} c_1, \text{ if } |\underline{c}| = 2r \text{ and } c_2 + c_3 \geq c_1 + c_4, \\ r - c_4, \text{ if } |\underline{c}| = 2r \text{ and } c_2 + c_3 \leq c_1 + c_4, \\ 0, \text{ otherwise} \end{cases}$

  [Fak09, Lemma 5.1]

- $\frac{1}{r} \deg \Phi_4^{\mathcal{G}}(\underline{c}) = \begin{cases} \min\{\frac{1}{r}c_1, 1 - \frac{1}{r}c_4\}, \text{ if } \frac{1}{r}|\underline{c}| = 2, \\ 0, \text{ otherwise} \end{cases}$

  [GG12, Lemma 2.6],[AS08, Lemma 2.2]

- $\deg \Phi_4^\Sigma(\underline{c}) = = \begin{cases} \min\{c_1, r - c_4\}, \text{ if } |\underline{c}| = 2r, \\ 0, \text{ otherwise} \end{cases}$

  [Fed11, Proposition 4.2],[BM10, EKZ10]

It is immediate that these three formulae coincide.